\documentclass[11pt]{amsart}
\usepackage{graphicx, amsfonts, amssymb}
\usepackage{mathrsfs, amsmath, amsthm}
\usepackage{verbatim}
\usepackage{relsize}
\usepackage{epsfig}
\usepackage{color, colortbl}
\usepackage{comment}
\usepackage{mymacros}
\usepackage{xfrac,bigints}
\usepackage{enumerate}
\usepackage{hyperref}
\hypersetup{
	colorlinks=true,
	linkcolor=blue,
	filecolor=magenta,      
	urlcolor=cyan,
	pdftitle={Optimal Domain for Tg Operator},
	pdfpagemode=FullScreen,
}

\theoremstyle{definition}

 \theoremstyle{remark}

\newtheorem{remark}{Remark}

\numberwithin{equation}{section}

\newcommand{\D}{\mathbb D}

\newcommand{\ma}{\mathbb}

\newcommand{\abs}[1]{\lvert#1\rvert}

\def\BMOA{\text{BMOA}}
\newcommand{\HOl}[1]{\text{Hol}\left( #1\right)}

\begin{document}

\title{Optimal Domain of generalized Volterra operators}
\author[C. Bellavita]{C. Bellavita$^1$ }
\author[V. Daskalogiannis]{V. Daskalogiannis$^{1,2}$}
\author[G. Nikolaidis]{G. Nikolaidis$^{1}$}
\author[G. Stylogiannis]{G. Stylogiannis$^{1}$}

\email{cbellavita@math.auth.gr}
\email{vdaskalo@math.auth.gr}
\email{nikolaidg@math.auth.gr}
\email{stylog@math.auth.gr}

\address{$^1$  Department of Mathematics, Aristotle University of Thessaloniki, 54124, Thessaloniki, Greece.}
\address{$^2$  Division of Science and Technology, American College of Thessaloniki, 17 V. Sevenidi St., 55535, Pylea, Greece.}

\thanks{This research project was supported by the
Hellenic Foundation for Research and Innovation (H.F.R.I.) under the '2nd
Call for H.F.R.I. Research Projects to support Faculty Members \& Researchers' (Project Number: 4662).}



\begin{abstract}
For $g$ in $\BMOA$, we consider the generalized Volterra operator $T_g$ acting on Hardy spaces $H^p$. This article aims to study the largest space of analytic functions, which is mapped by $T_g$ into the Hardy space $H^p$. We call this space the optimal domain of $T_g$ and we describe its structural properties. Motivation for this comes from the work of G. Curbera and W. Ricker \cite{Curbera2011} who studied the optimal domain of the classical Ces\'aro operator.
\end{abstract}
\maketitle
\section{Introduction}
 Let $\D$ be the unit disk of the complex plane  and $\Hol(\ma{D})$ be the space of analytic functions on $\D$. 
We consider the Hardy spaces of the unit disc
\[
 H^p=\{f\in \Hol(\ma{D}):\; \|f\|_p= \sup_{0\leq r<1}M_p(r,f)<+\infty\}\, ,
\]
where $1\leq  p<\infty$ and 
$M_p(r,f)=\left(\int_{0}^{2\pi}|f(re^{it})|^p\,\frac{dt}{2\pi}\right)^{1/p}$. For $p=\infty$, let $H^\infty$ denote the space of all bounded analytic functions on $\D$, with norm given by $\|f\|_\infty= \sup_{z \in \mathbb D}|f(z)|<+\infty$. See \cite{Duren1970} for more information on the Hardy spaces. The classical space $\BMOA=\text{BMO} \cap \Hol(\D)$, can be equivalently described as the {\it M\"{o}bius invariant} subspace of the Hardy spaces, i.e.
$$
\BMOA=\{g\in H^p:\; \|g\|_{\BMOA,p}=\sup_{a\in\ma{D}}\|g\circ \phi_a-g(a)\|_p<+\infty\}\, ,
$$
where $\phi_a(z)=\frac{a-z}{1-\overline{a}z},\;z\in\D$, are the M\"obius automorphisms of the disc. BMOA is a Banach space under the norm 
$$\|g\|_*=|g(0)|+\|g\|_{\BMOA,2}.$$
See \cite[Chapter VI]{garnett2006} and \cite{Girela2001} for basic properties of $\BMOA$, including the fact that the above equivalent definition does not depend on $p\in[1,+\infty).$

For $f,\,g \in \Hol(\ma{D})$, consider the generalized Volterra operator
$$
T_g(f)(z)= \int_{0}^{z}f(\zeta)g'(\zeta)\,d\zeta\ , \quad z \in \mathbb D\, .
$$

For $g(z)=z$, it reduces to the classical Volterra operator. For a general symbol $g$, Ch. Pommerenke \cite{Pommerenke1977} proved that $T_g\colon H^2\to H^2$ is bounded if and only if $g\in \BMOA$. Aleman and Siskakis \cite{Aleman1995} extended this result to all  $H^p$, for $1\leq p<+\infty$, by proving that boundedness of $T_g$ on $H^p$ happens for the same class of symbols $g$. Three years later, Aleman and Cima \cite{Aleman2001} completely characterized the functions $g$ for which the operator $T_g\colon H^p\to H^q$ acts boundedly, for all different choices of $0<p,\,q<+\infty$. We also consider its {\it companion} operator
$$
S_g(f)(z)=\int_{0}^{z}f'(\zeta)g(\zeta)\,d\zeta\ , \quad z \in \mathbb D\, ,
$$
and the multiplication operator $M_g(f)(z):= f(z)g(z)$. We say that a function $g$ is a multiplier for the space $X$, if $M_g:X\to X$ is bounded and we denote the set of multipliers of $X$ by $\mathcal{M}(X)$.  An integration by parts gives
$$
 T_g(f)(z)=M_g(f)(z)-g(0)f(0)-S_g(f)(z)\,.
$$ 
 It is clear that if two of the three operators are bounded, then the third one is also bounded. However, it is possible that operator $T_g$ is bounded even if $S_g$ and $M_g$ are not, due to cancellation.

The class of operators $T_g$ is a natural generalization of the well-known Ces\'aro operator $C$. Indeed, if $g(z)=-\log(1-z)$, then
$$
T_g(f)(z)= \int_{0}^{z}\frac{f(\zeta)}{1-\zeta}\,d\zeta=zC(f)(z)\, .
$$
It is known that the Ces\'aro operator is bounded on $H^p$, for all $0< p<\infty$, see \cite{Miao1992}, \cite{Siskakis1987}.
However, one can find functions $f$ not belonging to $ H^p$, for which $C(f)$ is still in $H^p$. In \cite{Curbera2011} and \cite{curbera2012}, Curbera and Ricker introduced the optimal domain for the Ces\'aro operator $[C,H^p]$, consisting of all analytic functions that $C$ maps into $H^p$. They studied the properties of the space $[C,H^p]$, and proved that it is a Banach space which strictly contains $H^p$. Furthermore, they showed that  $[C,H^{p_2}]\subsetneqq [C,H^{p_1}]$ for $1\leq p_1<p_2<+\infty$. In addition, they prove that the space of multipliers $\mathcal{M}([C,H^p])$ is $H^\infty$.

Following the notation used in \cite{Curbera2011}, we define the optimal domain of the generalized Volterra operator $T_g$ acting on $H^p$, as follows
    $$
    [T_g,H^p]:=\left\lbrace f \in \Hol(\D):\;  T_g(f)\in H^p\right\rbrace\,.
    $$
It is clear that $[T_g,H^p]$ is a linear subspace of $\Hol(\D)$. If $g$ is constant, then trivially $[T_g,H^p]=\Hol(\D)$. For $g$ non-constant, we equip the space with the norm
    $$
    \|f\|_{[T_g,H^p]}:=\|T_g(f)\|_{p}\ .
    $$ 
 In what follows, we study some basic properties of this normed space.

\section{Main results and plan of the paper}
From now on, we always assume that $g\in \BMOA$ and is not constant. 

\begin{thm}\label{Theorem Tg is a Banach Space}
Let $1\leq p<\infty$ and $g\in \BMOA$. Then $[T_g,H^p]$ equipped with the norm $\|\cdot\|_{[T_g,H^p]}$, is a Banach space.
\end{thm}
The optimal domain $[T_g,H^p]$ is a considerably large space. Clearly, for $g\in BMOA$, $H^p\subset [T_g,H^p]$. We show that this inclusion is strict.
\begin{thm}\label{Theorem H^p strict subset of optimal domain.}
    Let $1\leq p<\infty$ and $g\in \text{BMOA}$. Then $ H^p$ is strictly smaller than $[T_g, H^p]$. Moreover, for every $f \in H^p$
    $$
    \|f\|_{[T_g,H^p]}\leq C \|f\|_p\,  ,
    $$
where $C=\|T_g\|_{H^p\to H^p}$.
\end{thm}
This theorem is a direct consequence of the fact that the operators $T_g$ are never bounded from below in the Hardy spaces, see \cite{Anderson2011} and \cite{Panteris2022}.
Hence, a natural question arises: to what extent is the Hardy space $H^p$ smaller, compared to $[T_g,H^p]$? In Propositions \ref{Proposition 1/psi} and \ref{proposition g' carleson example}, we provide specific examples of functions $g\in \BMOA$ for which $[T_g,H^p]$ is strictly larger than $H^p$. In the same direction, we prove the following result; 
\begin{thm}\label{Theorem intersection }
Let $1\leq p<\infty$. Then
 $$
 H^p= \bigcap_{g\in \BMOA}[T_g,H^p]\, .
 $$    
\end{thm}
As we will see in the next section, there exists a connection between the optimal domain and the non-radial weighted Bergman spaces. For this reason, all the issues regarding non-radial weighted Bergman spaces also appear while studying $[T_g,H^2]$. For example, it is unclear for which functions $g$ the polynomials are dense in $[T_g,H^2]$.  
 
Since $[T_g,H^p]$ are Banach spaces of holomorphic functions, it makes sense to find out whether they share properties common to other spaces of analytic functions. We show that although $ H^p$ is strictly smaller than $[T_g, H^p]$, the space of pointwise multipliers does not change. 

\begin{thm}\label{Theorem multipliers introduction}
Let $1\leq p<\infty$ and $g\in \BMOA$. The multiplication operator $$M_h:\,[T_g,H^p]\to [T_g,H^p]$$ is bounded, if and only if $h \in H^\infty$. Moreover, there exists a constant $C>0$, such that $C\|h\|_{\infty}\leq \|M_h\|\leq \|h\|_{\infty}$.
\end{thm}

Finally, we prove that for a fixed symbol $g\in BMOA$, the optimal domains inherit the inverse inclusion property of the spaces $H^p$, with the inclusions being strict.
\begin{thm}\label{Theorem nikos2}
Let $g \in \BMOA$. For $1\leq p_1<p_2<+\infty$, we have that
$$
[T_g,H^{p_2}]\subsetneqq [T_g,H^{p_1}]\, .
$$
\end{thm}

The rest of the article is divided into three sections. In section \ref{Section 3}, we investigate the structural properties of the optimal domain. We prove Theorems \ref{Theorem Tg is a Banach Space}, \ref{Theorem H^p strict subset of optimal domain.}, and \ref{Theorem intersection }, we study the density of polynomials in $[T_g, H^2]$, and we discuss a connection between the optimal domain and the non-radial weighted Bergman spaces. Additionally, we provide specific examples of classes of symbols $g$, for which the optimal domain is significantly larger than $H^p$.
In section \ref{Multipliers} we turn our attention to the pointwise multipliers of the optimal domain, while we devote  section \ref{Inclusion section} to the, more technical, proof of Theorem \ref{Theorem nikos2}.

As usual, $A \sim B$ means that there are positive constants $c$ and $d$ such that $cB\leq A\leq dB$. Moreover, $p \in \mathcal{X}\setminus \mathcal{Y}$ means that the element $p$ belongs to the set $\mathcal{X}$ but not to $\mathcal{Y}$.

\section{The structure of the optimal domain}\label{Section 3}
 By definition, it is clear that $[T_g,H^p]$ is a linear subspace of $\Hol(\mathbb D)$. For this reason, our primary objective is to establish that $[T_g,H^p]$ is a Banach space. We recall from \cite{Anderson2011}, that $T_g$ is injective for non constant symbols in $\text{Hol}(\ma{D})$ and consequently $\|\cdot \|_{[T_g,H^p]}$ is a true norm for $p\geq 1$. 
\begin{proof}[Proof of Theorem \ref{Theorem Tg is a Banach Space}]
First of all, we prove that the point evaluations of $[T_g,H^p]$ are uniformly bounded in every compact subset $K$ of the unit disk. 

Let $z_0\in\D$. By \cite[p. 36]{Duren1970} and by applying Cauchy integral formula, we deduce that there exists a positive constant $C>0$ such that
    \begin{equation}\label{E:derivative function}
        |h'(z_0)|\leq \frac{C}{(1-|z_0|)^{1+\frac{1}{p}}}\|h\|_{p} \,
    \end{equation}
 for every $ h\in H^p$.
Let $f\in [T_g,H^p]$. Then $T_g(f)\in H^p$, so
  \[
    \begin{split}
      |f(z_0)\cdot g'(z_0)|&=\,|T_g'(f)(z_0)|
    \\
   &\leq \, \dfrac{C}{(1-|z_0|)^{1+\frac{1}{p}}}\,\|T_g(f)\|_{p}
   \\
   &=\,\dfrac{C}{(1-|z_0|)^{1+\frac{1}{p}}}\,\|f\|_{[T_g,H^p]}\, .
    \end{split}
    \]
Therefore, the boundedness of the point evaluations reduces to whether the above inequality can be divided  by $|g'(z_0)|$.

Following the lines of \cite[Section 5.4]{duren2004}, let $K$ be a compact subset of $\ma{D}.$ If $z_0\in K$ and $g'(z_0)\neq 0$, then
\[
|f(z_0)|\leq \frac{C}{|g'(z_0)|(1-|z_0|)^{1+\frac{1}{p}}}\,\|f\|_{[T_g,H^p]}\,.
\]
Thus by continuity of $g'$, the point evaluations are bounded in some neighbourhood of $z_0.$ If $g'(z_0)=0$, then there exists an annulus $A_{z_0}$ centered at $z_0$, such that     
$$
    \min_{\zeta\in A_{z_0}}|g'(\zeta)|>\varepsilon>0\, .
    $$ 
 Since $|f|$ is subharmonic, by applying the sub-mean value property, denoting by $dA(w)$ the normalized area measure on $\D$ and by $\abs{A_{z_0}}$ the area of the annulus $A(z_0)$, we have 
    \begin{align*}
        |f(z_0)|&\leq \frac{1}{|A_{z_0}|}\int_{A_{z_0}}|f(w)|\,dA(w)\\
        &\leq \frac{1}{|A_{z_0}|}\int_{A_{z_0}}\frac{C}{(1-|w|)^{1+\frac{1}{p}}\min_{\zeta\in A_{z_0}}|g'(\zeta)|}\,\|f\|_{[T_g,H^p]}\,dA(w)\\
        & =\frac{\|f\|_{[T_g,H^p]}}{\varepsilon|A_{z_0}|}\int_{A_{z_0}}\frac{1}{(1-|w|)^{1+\frac{1}{p}}}\,dA(w)= C\,\|f\|_{[T_g,H^p]}\, .
    \end{align*}
    where $C>0$ depends only on $z_0$. By continuity of $f$, the point evaluations are uniformly bounded also in a neighbourhood of $z_0$. Since $K$ is compact and there exists a finite cover of open sets. It follows that the point evaluations are uniformly bounded on $K$.
    
       To conclude, we prove the completeness of $\left( [T_g,H^p], \|\cdot\|_{[T_g.H^p]}\right)$. Let $\{f_n\}_n$ be a Cauchy sequence in $[T_g,H^p]$, that is
    \begin{equation}\label{equi2}
        \|f_n-f_m\|_{[T_g,H^p]}=\|T_g(f_n-f_m)\|_{p}\rightarrow 0\,,
    \end{equation}
as $m,n\to \infty$.
By the uniform boundedness of point evaluations in $[T_g,H^p]$ and (\ref{equi2}), $f_n$ is uniformly Cauchy in every compact subset of $\mathbb{D}\,$. Therefore, there exists a function $f \in \Hol(\mathbb D)$, such that
$f_n\to f$, locally uniformly on $\D$. Since $H^p$ is also a Banach space, (\ref{equi2}) implies that there exists a function $k\in H^p$ such that $T_g(f_n)$ converges in norm to $k$, and since the point evaluation for the derivative of the Hardy spaces are uniformly bounded on compact sets, $T_g'(f_n)(z)\rightarrow k'(z)$, for every $z \in \mathbb D$.

Combining the above ideas, for every $z \in \mathbb D$,
    \begin{align*}
        k'(z)&=\lim_{n\to+\infty} f_n(z)\cdot g'(z) =g'(z)\cdot f(z)
    \end{align*}
    and, by taking the primitives,  $k=T_g(f)$ from which the statement of the theorem follows.
\end{proof}
It is clear that if $g \in \BMOA$, then $H^p\subset [T_g,H^p]$. We prove that the inclusion is strict. In this direction, we use the fact that $T_g$ is never bounded from below in $H^p$, extending a well known result from \cite{Anderson2011} to the case $1\leq p<+\infty$.

\begin{lem}\label{bounded below Tg Lemma}
    Let $1\leq p<+\infty$ and $g\in \BMOA$. Then $T_g\colon H^p\to H^p$ is not bounded below. 
\end{lem}
\begin{proof}
Assume that $T_g$ is bounded from below on $H^p$. This means that there exists a positive constant $C>0$ such that
    $$
    \|T_g(f)\|_p\geq C\|f\|_p\qquad \forall f\in H^p \, .
    $$
We first consider $1\leq p\leq 2$. Our assumption implies that, for $f(z)=z^n$ for which $\|z^n\|_{p}=1$, there exists a constant $C>0$ such that
    \begin{equation}\label{equi1 Lemma 3.1}
        \|T_g(z^n)\|_p\geq C\|z^n\|_p=C.
    \end{equation}
However, taking into account the inequality

    $$\|T_g(f)\|_p\leq \|T_g(f)\|_2,$$
for every $f \in H^2$, and, \cite{Anderson2011}, $\|T_g(z^n)\|_2\to0$, we have that $\|T_g(z^n)\|_p\to 0$, which contradicts \eqref{equi1 Lemma 3.1}.

For the case $p>2$, by an application of the Cauchy-Schwarz inequality, we have that
\begin{align*}
 \|T_g(z^n)\|_p^p\,&\leq \, \|T_g(z^n)\|_{2p-2}^{p-1}\,
 \|T_g(z^n)\|_2 \\
 &\leq\, C\|z_n\|_{2p-2}^{p-1}\,\|T_g(z^n)\|_2 \,=\, C\|T_g(z^n)\|_2\,,
\end{align*}
where in the second inequality we have used that $T_g$ is bounded in $H^{2p-2}$. Once more, using the fact that $\|T_g(z^n)\|_2$ converges to zero, we conclude that $T_g$ cannot be bounded from below.
    \end{proof}

Using this result, we prove that $H^p$ is a proper subset of $[T_g, H^p]$.
\begin{proof}[Proof of Theorem \ref{Theorem H^p strict subset of optimal domain.}]
Since for $g \in \BMOA$, $T_g$ is bounded \cite{Aleman1995}, it follows that $ H^p \subset [T_g, H^p]$ and, in particular, that for every $f \in H^p$
$$
\|f\|_{[T_g,H^p]}=\|T_g(f)\|_{p}\leq \|T_g\|_{H^p \to H^p}\|f\|_{p}\, .
$$
We assume that there exists a function $g \in \BMOA$ for which $[T_g, H^p]= H^p$.
In that case, the identity operator $I_d : H^p\to [T_g,H^p]$ is a continuous, onto operator between two Banach spaces and by the open mapping theorem
$$
\|f\|_p\sim \|f\|_{[T_g,H^p]}\, .
$$
In particular, for every $f \in H^p$, there exists a positive constant $C$, such that
\begin{equation}\label{E: bounded from below}
\|T_g(f)\|_p\geq C\|f\|_p\, .    
\end{equation}
However, this is a contradiction, since $T_g$ is not bounded from below in $H^p$, as we showed above.
\end{proof}

\subsection*{Examples}For some particular choices of $g \in \BMOA$, we can explicitly exhibit  functions in $[T_g, H^p]\setminus H^p$. To provide our first example, we consider
$$
H^p(\omega)=\left\lbrace f \in \text{Hol}(\mathbb D)\ \text{ such that }\ f \omega=g \in H^p(\D) \right\rbrace 
$$
 where $\omega$ is an outer function. From \cite[Proposition 2.1]{Curbera2011} we know that
\begin{equation}\label{subset H^p}
 H^p(\omega)\subset H^p(\D) \iff  1/\omega\in H^\infty(\D)\, .   
\end{equation}
Let $g$ be an element of $\BMOA_{log}\subset \BMOA$. This weighted BMOA space appears in the study of the generalized Volterra operator. In specific, the set of symbols $g\in \Hol(\ma{D})$ for which $T_g$ is bounded on $\BMOA$, is exactly $\BMOA_{log}$ \cite{Siskakis1999}. We remind the reader that a function $f$ is in $\BMOA_{log}$ if
\[
\sup_{I\subseteq \partial\D }\, \frac{\left(\log\dfrac{e}{|I|}\right)^2}{|I|}\,\int_{S(I)} \vert f'(z)\vert^2(1-|z|^2)\,dA(z)\,<\,\infty\,,
\]
where $I$ is an arc on the unit circle with length $\abs{I}$, and $S(I)$ is the usual {\it Carleson box}, i.e.  $S(I)= \{z\in\D: \frac{z}{|z|}\in I, \, \text{and} \, 1-|I|<|z|<1\}$.

\begin{prop}\label{Proposition 1/psi}
Let $g \in \BMOA_{log}$ and $\psi$ be an outer function. If $1/\psi \in \text{BMOA}\setminus H^{\infty}$, then there exists $f \in H^p(\psi)$ such that  $f \in [T_g,H^p]\setminus H^p$.
\end{prop}
\begin{proof}
If $f \in H^p(\psi)$, there exists $h \in H^p(\D)$ such that $f=h/\psi$. In particular
$$
T_g(f)(z)= \int_{0}^z g'(\xi) f(\xi) d\xi = \int_{0}^z g'(\xi)\frac{h(\xi)}{\psi(\xi )} d\xi=T_k(h)\, $$
where
$$
k(z)=\int_0^z \frac{g'(\xi)}{\psi(\xi)}d\xi=T_g\left(\frac{1}{\psi}\right)(z)\, .
$$
From \cite[Theorem 3.1]{Siskakis1999}, $k\in \BMOA$ and so  $T_k(h)\in H^p$. Therefore $T_g$ maps $H^p(\psi)$ into $H^p$,
that is $H^p(\psi) \subset [T_g, H^p]$. To finish the proof, notice that after \eqref{subset H^p}, for this $\psi$ we can always find a function $f\in H^p(\psi)\setminus H^p$.

We note that an example of a pair of functions $(\psi,\,f)$ as described in the statement of this proposition, is given by
$$
\psi(z)=\frac{1}{3-\log(1-z)}\,,\quad  f(z)=(1-z)^{-1/p} \left( \log\left(\frac{1}{1-z}\right)\right)^{1-1/2p}\, .
$$
\end{proof}
As an alternative example, we show that it is possible to identify classes of functions larger than $H^p$, yet still contained within $[T_g, H^p]$, by considering symbols $g \in \BMOA$, as described below.
\begin{prop}\label{proposition g' carleson example}
    Let $I\subset \ma{T}$ be an arc centered at $a\in \ma{T}$. We assume that $g\in \BMOA$ is such that in the Carleson square $S(I)$
    $$
    \sup_{z \in S(I)}\abs{g'(z)}<C
    $$
for some positive constant $C$.
Then $H^p\subsetneq H^p(z-a) \subset  [T_g,H^p]$.
\end{prop}
\begin{proof}
It suffices to show that
$$
k(z)=\int_0^z \frac{g'(\xi)}{\xi-a}d\xi \in \BMOA\, 
$$
and, by following the steps of the proof of Proposition \ref{Proposition 1/psi}, we prove the right hand side of the desired inclusion. In order to do so, we use  a description  of $\BMOA$ in terms of Carleson measures, see \cite[Theorem VI.3.2]{garnett2006}. For simplicity we assume that $|I|\leq 2\delta$ and, for a fixed $J \subset \ma{T}$, we write 
$A=S(J)\cap S(I)$ and $B=S(J)\cap (\ma{D}\setminus S(I))$.
Then
{\small\[
\begin{split}
    \norm{k}&^2_{\text{BMOA}}\sim \sup_{J} \frac{1}{\abs{J}} \int_{S(J)}\left|\frac{1}{\xi-a}\right|^2\abs{g'(\xi)}^2(1-\abs{\xi}^2)dA(\xi)\\
    &=\sup_{J} \frac{1}{\abs{J}}  \left( \int_{A} +\int_{B}\right) \bigg|\frac{1}{\xi-a}\bigg|^2\abs{g'(\xi)}^2(1-\abs{\xi}^2)dA(\xi)\\
    &\leq \sup_{J} \frac{1}{\abs{J}}  \left( C^2 \int_{A}\bigg|\frac{1}{\xi-a}\bigg|^2(1-\abs{\xi}^2)dA(\xi) + \frac{1}{\delta^2} \int_{B} \abs{g'(\xi)}^2(1-\abs{\xi}^2)dA(\xi)\right)\\
    &\leq C^2 \sup_{J} \frac{1}{\abs{J}} \left(\int_{S(J)}\bigg|\frac{1}{\xi-a}\bigg|^2 (1-\abs{\xi}^2)dA(\xi) +\frac{1}{\delta^2} \int_{S(J)} \abs{g'(\xi)}^2(1-\abs{\xi}^2)dA(\xi)\right)\\
    &\leq C' \norm{\log(\cdot-a)}_{\text{BMOA}}^2 + C'' \norm{g}_{\text{BMOA}}^2\,.
    \end{split}
\]}
The proof is complete once one notices that $(z-a) \in H^\infty$, while $1/(z-a) \notin H^\infty$ and, consequently, $H^p\subsetneq H^p(z-a))$, due to \eqref{subset H^p}.

A trivial example of a function $g$ satisfying the requirement of the proposition is given by $g(z)=-\log(1-z)$, with $S(I)$ centered at $a=-1$.
\end{proof}
Next, we turn our attention to Theorem \ref{Theorem intersection }. 

\begin{proof}[Proof of Theorem \ref{Theorem intersection }]
    
    Since $T_g(H^p)\subset H^p$ for each $g\in \BMOA$, we have that
    $$
    H^p\subseteq \bigcap_{g\in \BMOA}[T_g,H^p]\, .
$$
    

Now we verify the inverse inclusion. Assume $\displaystyle{f\in \bigcap_{g\in \BMOA}[T_g,H^p]}$ and consider the companion operator $S_f(g)(z)=\int_{0}^{z}g'(\zeta)f(\zeta)d\zeta$ with $g\in \BMOA$. The assumption on $f$ and the relation
 $$
 T_g(f)(z)=S_f(g)(z)=\int_{0}^{z}g'(\xi)f(\xi) d\xi
 $$
 yield, $S_f(g)\in H^p$ for every $g\in \BMOA$. By the closed graph theorem, this implies that $ S_f\colon \BMOA\rightarrow H^p$ is bounded, and we show that this implies $f\in H^p$.
Let $f(z)=\sum_{k= 0}^{\infty}a_k z^k$ be the Taylor series for $f$ and $f_n(z)=z^n$, $n\geq 0$ the monomials. Then
$$
S_f(f_n)(z)=\int_{0}^z f(\zeta)(\zeta^n)'d\zeta=\sum_{k=0}^{\infty}a_k\frac{n}{n+k}z^{n+k},
$$
and consider the functions
$$F_n(z)=\frac{1}{z^n}S_f(f_n)(z)=\sum_{k=0}^{\infty}a_k\frac{n}{n+k}z^{k},\quad n\geq 0.$$
We claim that 
\begin{equation}\label{equioffinSf(zn)}
    f(z)=\lim_{n\to\infty}F_n(z) \qquad z \in \mathbb D.
\end{equation}

Indeed, $\displaystyle{f(z)-F_n(z)=\sum_{k=0}^{\infty}\frac{k a_k}{n+k}z^k}$, and for a fixed $z$ and $\varepsilon>0$, there is $M\in\ma{N}$ such that $$\sum_{k=M+1}^{\infty}|a_k||z|^k<\varepsilon.$$
Also there is $M\in\ma{N}$, such that
$$\left|\sum_{k=0}^M\frac{ka_k}{n+k}z^k\right|<\varepsilon\quad \text{ for }n\geq M.$$
Thus,
\begin{align*}
    |f(z)-F_n(z)|&\leq \left|\sum_{k=0}^M\frac{ka_k}{n+k}z^k\right|+\sum_{k=M+1}\frac{k|a_k|}{n+k}|z|^k\\
    & \leq  \left|\sum_{k=0}^M\frac{ka_k}{n+k}z^k\right|+\sum_{k=M+1}|a_k||z|^k\\
    & < 2\varepsilon
\end{align*}
for $n\geq M$ and \eqref{equioffinSf(zn)} follows. Next, since $S_f\colon \BMOA\to H^p$ is bounded and $\|f_n\|_*\leq 1$ for each $n$, we have
\begin{align*}
    \|F_n\|_p&=\|z^nF_n\|_p=\|S_f(f_n)\|_p\\
    & \leq \|S_f\|_{\BMOA\to H^p}\|f_n\|_*\\
    & \leq \|S_f\|_{\BMOA\to H^p}.
\end{align*}
From Fatou's Lemma then, for $r<1$, we have
\begin{align*}
M_p^p(r,f)&=\int_{0}^{2\pi}\liminf_{n\to+\infty}|F_n(re^{i\theta})|^p\frac{d\theta}{2\pi}\leq \liminf_{n\to+\infty}\int_{0}^{2\pi}|F_n(re^{i\theta})|^p\frac{d\theta}{2\pi}\\
& \leq \|F_n\|_p^p\leq \|S_f\|_{\BMOA\to H^p}^p
\end{align*}
and $\|f\|_p^p=\lim_{r\to 1}M_p^p(r,f)\leq \|S_f\|_p^p,$ thus $f\in H^p$ and the proof is finished.
\end{proof}
 
\begin{remark}
    We note that the above proof implies that  $S_f\colon \BMOA\rightarrow H^p$ is bounded if and only if $f\in H^p$. Moreover,
$$    \|f\|_{p}\sim ||S_f||_{\BMOA\to H^p}\, .
    $$
\end{remark}

\subsection{\texorpdfstring{$[T_g,H^2]$ }\empty as a Bergman space}In the case $p=2$, we can identify the space $[T_g,H^2]$ with a non-radial weighted Bergman space. Let $w\colon \ma{D}\rightarrow [0,+\infty)$ be a weight function. The weighted Bergman space $\cA^2_w$ is defined as 
$$
\cA^2_w=\left\lbrace f \in \text{Hol}(\mathbb D):\;  \|f\|_{\cA^2_w}^2:=\int_{\mathbb D}\abs{f(z)}^2w(z)dA(z)<\infty\right\rbrace\, ,
$$
where $dA(z)$ is the Lebesgue area measure of $\D$.
For the standard weights $w(z)=(\alpha+1)(1-|z|^2)^\alpha$,  the weighted Bergman space $\cA^2_w$ is denoted as $A_\alpha^2.$ For a complete description of these spaces, we refer the reader to \cite{duren2004} and \cite{hedenmalm2000}.
 Let $g \in \BMOA$ and consider the weight function $w(z)=\abs{g'(z)}^2\left( 1-\abs{z}^2\right)$. Then
$$
[T_g,H^2]=\cA^2_{w}\, .
$$
Indeed, by the {\it Littlewood - Paley} formula \cite[Lemma 3.2]{garnett2006}, if $f \in [T_g,H^2]$, then
$$
\|f\|_{[T_g,H^2]}=\|T_g(f)\|_{2}\sim \int_{\mathbb{D}}\abs{f(\xi)}^2\abs{g'(\xi)}^2\left( 1-\abs{\xi}^2\right)dA(\xi)\,.$$

Curbera and Ricker \cite{Curbera2011}, proved that the polynomials are dense in the optimal domain $[C,H^2]$. However, things are more complicated when considering a general symbol $g$. It is easy to see that the density of the polynomials in $[T_g,H^2]$ is equivalent to the density of $H^2$ in $[T_g,H^2]$.

\begin{prop}
Let $g \in \BMOA$. The polynomials are dense in $[T_g,H^2]$ if and only if $H^2$ is dense in $[T_g,H^2]$.    
\end{prop}
\begin{proof}
    It is clear that if the polynomials are dense in $[T_g,H^2]$, then also $H^2$ is dense in $[T_g,H^2]$. On the other hand, let us suppose that $H^2$ is dense in $[T_g,H^2]$. For given $\varepsilon>0$, $f\in [T_g,H^2]$, there exists $h \in H^2$ such that $\|f-h\|_{[T_g,H^2]}<\varepsilon/2$. As the polynomial are dense in $H^2$, see for instance \cite[Theorem 3.3]{Duren1970}, we can find a polynomial $p$ such that $\| h-p\|_{2}<\varepsilon/2$. Therefore
    \begin{align*}
        \|f-p\|_{[T_g,H^2]}&\leq \|f-h\|_{[T_g,H^2]} +\|h-p\|_{[T_g,H^2]}\\
        & \leq \|f-h\|_{[T_g,H^2]} + C\|h-p\|_{2}\\
        &\leq(1+ C)\varepsilon\ ,
    \end{align*}
    where in the second inequality, we use the fact that $T_g$ is bounded on $H^2$ and $C=C(g)$.
\end{proof}

In general, very little is known about the non-radial weighted Bergman spaces, see for example \cite{duren2004} and \cite{Aleman2009}. Nevertheless, using the theory of the standard weighted Bergman spaces $A^2_\alpha$, we are able to provide a specific example of a function $g \in \BMOA$, for which the polynomials are not dense in $[T_g,H^2]$. The proof follows the ideas of \cite[p. 138]{duren2004}.

\begin{prop}\label{polynomial density in Tg,H^2}
    Let $g(z)=\displaystyle{\int_0^z\exp\left( \frac{\xi+1}{\xi-1}\right)d\xi}$. The polynomials are not dense in $[T_g,H^2]$. 
\end{prop}
\begin{proof}
    Firstly, we recall from \cite[Theorem N]{Shapiro}, that the singular inner function
    $$S(z)=\exp\left( -\int_{\mathbb T}\frac{\xi+z}{\xi-z}\delta_1(\xi)\right)=\exp\left( \frac{z+1}{z-1}\right)$$
    is not cyclic in $A^2_1.$ Now, let
    $$g(z):=\int_{0}^{z}S(\xi)\,d\xi.$$
    It is true that $g$ belongs to $\BMOA$, see for instance \cite[Theorem 3.11]{Duren1970}. We prove that $1/S \in [T_g,H^2].$ Indeed, it is analytic in $\ma{D}$ and
    $$
\norm{1/S}_{[T_g,H^2]}^2\sim \int_{\mathbb D}(1-\abs{z}^2)dA(z)<\infty\, .
$$
On the other hand, $1/S$ cannot be approximated in $[T_g,H^2]$ by polynomials. If it was, then for every $\varepsilon>0$ there would exist a polynomial $p$ such that
{\small
\[
   \varepsilon> \norm{1/S- p}_{[T_g,H^2]}^2\sim \int_{\mathbb D}\abs{1-p(\xi)S(\xi)}^2(1-\abs{\xi}^2)dA(\xi)= \norm{1-pS}_{\cA^2_1}^2\, .
\]}
Therefore $S$ would be cyclic in $A^2_1$, which is a contradiction.
\end{proof}

\section{On pointwise multipliers and conformal invariance}\label{Multipliers}

 As we have proved in Theorem \ref{Theorem Tg is a Banach Space}, $\left([T_g,H^p], \|\cdot\|_{[T_g,H^p]}\right)$ is a Banach space. A next natural step is to describe its pointwise multipliers.
\begin{proof}[Proof of Theorem \ref{Theorem multipliers introduction}]         
Let $f\in [T_g,H^p]$ and $h \in H^\infty$. We note that
    $$
   T_g(hf) = \int_{0}^{z}h(\zeta)f(\zeta)g'(\zeta)d\zeta= S_h(T_g(f))
    $$
     Hence,
    \begin{align*}
   \|M_h(f)\|_{[T_g,H^p]}&=\|T_g(hf)\|_{p}=\|S_h(T_g(f))\|_{p}\\
   &\leq C\,\|T_g(f)\|_{p}\,=\, C\,\|f\|_{[T_g,H^p]}\,<\,\infty\,,
    \end{align*}
    where $C\,=\,\| S_h\|_{H^p \to H^p}\,= \|h\|_{\infty}$.

On the other hand, let us consider $h\in \text{Hol}(\ma{D})$, such that $ M_h\colon [T_g,H^p]\rightarrow [T_g,H^p]$
is bounded. Let $z\in\ma{D}$ and denote by $\lambda_z$ the point evaluation functional in $[T_g,H^p]$. As Theorem \ref{Theorem Tg is a Banach Space} shows, $\lambda_z$ is a bounded linear functional and since $[T_g,H^p]$ contains the constant functions, we have $\|\lambda_z\|\neq 0.$ For $f\in[T_g,H^p]$ with $\|f\|_{[T_g,H^p]}=1$ we have
$$\|\lambda_z\|:=\|\lambda_z\|_{[T_g,H^p]}\neq 0\qquad \forall z\in\ma{D}.$$
For $f\in [T_g,H^p]$, with $\|f\|_{[T_g,H^p]}=1$, we have that
\begin{align*}
    |h(z)f(z)|&=|\lambda_z(M_h(f))|\leq \|\lambda_z\|\|M_h\|.
\end{align*}
Thus, taking the supremum on $f\in [T_g,H^p]$ with $\|f\|_{[T_g,H^p]}=1$, we have
$$
|h(z)|\|\lambda_z\|\leq \sup_{\|f\|_{[T_g,H^p]}=1}|h(z)f(z)| \leq \|M_h\|\|\lambda_z\|
$$
which gives $h\in H^{\infty}$ since $\|\lambda_z\|\neq 0$.
\end{proof}

Curbera and Ricker also consider composition operators induced by automorphisms. They prove that in $[C,H^p]$ with $p\geq 1$, automorphisms do not induce bounded composition operators, see \cite[Proposition 3.6]{Curbera2011}. In the more general setting of $[T_g, H^p]$ and for specific symbols $g\in \BMOA$, it is possible to have bounded composition operators in the space $[T_g, H^p]$ (e.g. when $0<\inf_{z\in \mathbb{D}}\abs{g'(z)}\leq \sup_{z\in \mathbb{D}}\abs{g'(z)}<\infty$). However, even in these special cases it's not difficult to imagine that the space $[T_g,H^p]$ is not conformally invariant.  We use the result of \cite{Rubel1979} which asserts that the Bloch space $\cB$ is the largest conformally invariant space under certain natural assumptions. We recall that $\cB$ is defined as
$$
\cB=\left\lbrace f \in \Hol(\mathbb D)\colon \|f\|_{\cB}=\sup_{z \in \mathbb D} |f'(z)|(1-|z|^2)<\infty\right\rbrace\ . 
$$
For more information on $\cB$, we refer to \cite{zhu2007}.

\begin{cor}
  Let $1\leq p<\infty$ and $g\in \BMOA$.    The space $[T_g,H^p]$ is never conformally invariant.  
\end{cor}
\begin{proof}
Let us suppose that $[T_g,H^p]$ is conformally invariant. As point evaluation functionals in $[T_g,H^p]$ are uniformly bounded in compact subsets of $\ma{D}$, due to Theorem \ref{Theorem Tg is a Banach Space}, $[T_g,H^p]$ is a space of analytic functions satisfying the condition of \cite[Theorem]{Rubel1979}, that is there exists a linear functional $L$ on $[T_g,H^p]$ such that
$$|L(f)|\leq M\sup_{z\in K}|f(z)|$$
for all $f\in [T_g,H^p]$, for some $M>0$ and some compact set $K\subset \ma{D}$.
In this case, the identity operator $I_d\colon [T_g,H^p]\rightarrow \mathcal B$ is bounded. However this is never true, since $[T_g,H^p]$ always contains $H^p$ and there are functions which belong in $H^p\setminus \mathcal{B}$, which is a contradiction. For example, one such function is $f(z)=1/(1-z)^{2/3p}$.
\end{proof}

\section{Inclusion property}\label{Inclusion section}

 In this final section, we verify the {\it inverse inclusion} property. Since $H^{p_2} \subsetneq H^{p_1}$ for $p_1<p_2$, it is clear that $[T_g,H^{p_2}]\subset [T_g,H^{p_1}]$.
We will  prove that this inclusion is  strict for all $g \in \BMOA$.

It turns out that the zeros of $g'$ play an important role in this direction. It is known that, in general, $g'$ does not belong to the Nevanlinna class, see for instance \cite{Cohn1999} and \cite{Dyakonov2012}.
To provide an insight of the importance of the zeros, we first consider $g'=BG$, with $G$ non vanishing in $\mathbb{D}$ and $B$ a possibly infinite  Blaschke product, even if this case does not cover the general situation.

\begin{proof}[Proof of Theorem \ref{Theorem nikos2} when $g'=BG$]
 We assume that $[T_g,H^{p_2}]=[T_g,H^{p_1}]$. 

\noindent 
First of all, we consider the case $B \neq 0,1$.   For every $h \in H^{p_1}$, we consider 
$$
\phi=\frac{B}{g'}\, h' \in \Hol(\mathbb D)\ .
$$
Then
\begin{align*}
T_g(\phi)(z)=\int_{0}^{z}h'(t)B(t)dt=S_B(h)(z)\in H^{p_1}  
\end{align*}
since $B \in H^\infty$. Consequently, $\phi \in [T_g,H^{p_1}]$ and, because of our assumption, it is also in $[T_g,H^{p_2}]$. Therefore, for every $h \in H^{p_1}$,
$$
S_B(h)(z)=\int_{0}^{z}h'(t)B(t)\, dt=\int_{0}^{z}\frac{h'(t)B(t)}{g'(t)}g'(t)\, dt=T_g(\phi)(z) \in H^{p_2}\,,
$$
and by the closed graph theorem, the operator $S_B\colon H^{p_1}\to H^{p_2}$ is bounded. However, this is not possible. Indeed, by \cite[Theorem 2.2]{Anderson2011}, if $S_B\colon H^{p_1}\to H^{p_2}$ is bounded, then, for every $z \in \mathbb D$, 
$$
|B(z)|\leq C\frac{1/(1-|z|)^{1+1/p_{2}}}{1/(1-|z|)^{1+1/p_{1}}}=C(1-|z|)^{1/p_1-1/p_2}\ .
$$
In particular, the boundary value of $B$ would be zero, that is, $B$ would be the zero function, which contradicts our initial assumption.  

\vspace{11 pt}

On the other hand, if $B=1$ or equivalently $g'\neq 0$ on $\ma{D}$, we consider 
$$
\phi(z)=\frac{1}{(1-z)^{1/p_2+1}}\frac{1}{g'(z)} \in \Hol(\mathbb D)\, .
$$
Clearly, $\phi \in [T_g,H^{p_1}]\setminus [T_g,H^{p_2}]$, since
{\small\[
T_g(\phi)(z):=\int_0^z \frac{1}{(1-t)^{1/p_2+1}}dt=-p_2\left(\frac{1}{(1-z)^{1/p_2}}-1\right) \in H^{p_1}\setminus H^{p_2}\, .
\]}
\end{proof}

In order to face the general case, we recall the definitions of the Lipschitz and the Korenblum spaces, which are involved in our proof of Theorem \ref{Theorem nikos2}.
 If $0<\alpha\leq 1$, the Lipschitz space $\Lambda_\alpha$ is defined as
 $$
 \Lambda_\alpha=\left\lbrace f \in\Hol(\mathbb D):\; \|f\|_{\Lambda_\alpha}=\sup_{z \in \mathbb D}(1-|z|^2)^{1-\alpha} |f'(z)|<\infty\right\rbrace\,  .
 $$
 The Lipschitz space $\Lambda_\alpha$ plays a crucial role in the study of the boundedness of $T_g$ in the Hardy spaces. Indeed, in \cite[Theorem 1]{Aleman2001}, Aleman and Cima proved that, if $p_1<p_2$, then
\begin{equation}\label{E:AC}
    T_g: H^{p_1}\to H^{p_{2}} \; \text{is bounded} \iff g\in \Lambda_{\frac{1}{p_1}-\frac{1}{p_2}}\,  . 
\end{equation}
For an in-depth description of the spaces $\Lambda_\alpha$, we refer to \cite{zhu1993}.

For $0\leq \alpha<1$, we also consider the Korenbulm space $K_\alpha$, defined as
 $$
 K_\alpha=\left\lbrace f \in\Hol(\mathbb D) :\; \|f\|_{K_\alpha}=\sup_{z \in \mathbb D}(1-|z|^2)^\alpha |f(z)|<\infty\right\rbrace\, .
 $$
 We recall here some of the properties of $K_\alpha$ that we will use. First of all, it is clear that $f \in \Lambda_\alpha$ if and only if $f' \in K_{1-\alpha}$. Moreover, if $f \in K_{\alpha}$, then $f \in K_\beta$, for every $\beta>\alpha$.  For an introduction to Korenblum spaces, we refer the reader to \cite{korenblum1975}. We utilize a result which asserts that it is always possible to find a pair of holomorphic functions whose sum reaches the maximal possible growth, see \cite{Abakumov2013} and \cite{Xiao}.  More specific, let $\alpha\geq 0$. Then, there exist two holomorphic functions $f_1,f_2$ such that
\begin{equation}\label{E:Xiao}
|f_1(z)|+|f_2(z)|\sim (1-|z|)^{-a}
\end{equation}
for every $z \in \mathbb D$. Note that this implies that both $f_1$ and $f_2$ belong to $K_a$. 
 Using this result, we are able to characterize the multipliers mapping $K_\gamma$ to $K_\delta$. 

\begin{prop}\label{P:multipliers korenbul}
Let $0\leq \gamma<\delta$. Then 
$M_g:\;K_\gamma \to K_\delta$ is bounded, if and only if $\,g \in K_{\delta-\gamma}$.
\end{prop}
\begin{proof}
Let $g\in K_{\delta-\gamma}$. Then, for every $f \in K_\gamma$
\[
(1-|z|)^{\delta}|M_g(f)(z)|
\,=\, (1-|z|)^{\delta-\gamma}|g(z)|\; (1-|z|)^{\gamma}|f(z)|,\;\;z\in \D\,,
\]
hence $\norm{M_g(f)}_{K_\delta}<\infty$.

For the converse implication, let $M_g$ be bounded from $K_\gamma$ to $K_\delta$. After \eqref{E:Xiao}, we consider $f_1,\, f_2$  and $\alpha=\gamma$, so that $f_1,\,f_2\in K_\gamma$. By our assumption, $M_g(f_1),\,M_g(f_2)\in K_\delta$. Noticing that
$$
(1-|z|)^{\delta-\gamma}|g(z)|\,\leq\, C (1-|z|)^{\delta}\,(|f_1(z)|+|f_2(z)|)\,|g(z)|\,, 
$$
we conclude that $g \in K_{\delta-\gamma}$.
\end{proof}
  Before we proceed to the proof of Theorem \ref{Theorem nikos2}, we need to prove the following auxiliary result, which is crucial to our proof.
\begin{lem}\label{L:missing}
Let $g \in K_\gamma$ with $0<\gamma<1$ and let $0<\varepsilon<\min\{\gamma,1-\gamma\}$. Then, there exists a function $\phi \in \Hol(\mathbb D)$ such that $\phi \cdot g \in K_{\gamma+\varepsilon}\setminus K_{\gamma}$. 
\end{lem}
\begin{proof}
Let $\gamma'=\inf\{\beta\in(0,1)\colon g\in K_\beta\}$ ($\gamma'\leq \gamma$), and consider the Korenblum space $K_{\gamma'}$.
For $\alpha=\gamma-\gamma'$, we use the pair $f_1,\,f_2$ as in \eqref{E:Xiao}. Then $f_1\cdot g,\,f_2\cdot g\in K_\delta$ for every $1>\delta>\gamma$. Indeed, for $i=1,\,2$
\[
    |g(z)|\,|f_i(z)|\,(1-|z|)^\delta\,\leq\, C\,|g(z)|(1-|z|)^{\gamma'+(\delta-\gamma)}\,
\]
which is bounded due to the definition of $\gamma'$, since $\delta-\gamma>0$.
If $f_1\cdot  g \notin K_\gamma$, then $f_1$ is the desired function $\phi$. We have, therefore, reduced the problem to $f_1\cdot g\in K_\gamma$.

For the chosen $\varepsilon$, at least one of the functions $f_1\cdot g$ and $f_2\cdot g$ does not belong to $K_{\gamma-\varepsilon}$. To see this, let us assume on the contrary, that both $f_1\cdot g,\,f_2\cdot g \in K_{\gamma-\varepsilon}$.  Then,  noticing that
\begin{align*}
     |g(z)|(1-|z|)^{\gamma'-\varepsilon}\,&=\, |g(z)|(1-|z|)^{\gamma-\varepsilon}\,(1-|z|)^{-(\gamma-\gamma')}\\
     &=\, |g(z)|(1-|z|)^{\gamma-\varepsilon}\,(1-|z|)^{-\alpha}
\end{align*}
and using \eqref{E:Xiao} once more, we get that
\[
\begin{split}
\norm{g}_{K_{\gamma'-\varepsilon}}\,&\leq\, C\, \sup_{z \in \mathbb D} \,|g(z)|\,(1-|z|)^{\gamma-\varepsilon}\,\left( |f_1(z)|+|f_2(z)|\right)
\\
&\leq\,C\, \norm{f_1\,g}_{K_{\gamma-\varepsilon}} \,+\,\norm{f_2\,g}_{K_{\gamma-\varepsilon}} \,<\,\infty\,,
\end{split}
\]
which contradicts the definition of $\gamma'$.
 Without loss of generality, from now on we assume that $f_1\cdot g\notin K_{\gamma-\varepsilon}$.

For $\psi \in K_\varepsilon$ we consider the function 
$\psi\cdot f_1$. Notice first that

\[
\abs{\psi(z)\,f_1(z)\,g(z)}(1-|z|)^{\gamma+\varepsilon}=\abs{\psi(z)}(1-|z|)^\varepsilon\,\abs{f_1(z)\,g(z)}(1-|z|)^\gamma\, .
\]
Hence
\[
\norm{\psi\,f_1\,g}_{K_{\gamma + \varepsilon}}\leq \norm{\psi}_{K_\varepsilon} \,\norm{f_1\,g}_{K_\gamma}\,<\,\infty\,.
\]
We claim that we can always find such a function $\psi$, for which $\psi\cdot f_1\cdot g \notin K_\gamma$.

Indeed, if such $\psi$ did not exist, then for every function  $h\in K_\varepsilon$ we should have $h\cdot f_1\cdot g \in K_\gamma$ and consequently that $M_{f_1\cdot g}\colon K_\varepsilon\to K_\gamma$ is bounded. By Proposition \ref{P:multipliers korenbul},  this is equivalent to  $f_1\cdot g\in K_{\gamma-\varepsilon}$ which contradicts our  assumption. Hence, the desired function $\phi$ of the statement is given by $\phi=\psi\cdot f_1$.
\end{proof}

We are now ready to prove Theorem \ref{Theorem nikos2}.

\begin{proof}[Proof of Theorem \ref{Theorem nikos2}]
Let us assume that $[T_g,H^{p_2}]=[T_g,H^{p_1}]$, aiming to arrive at a contradiction.

Let $Z(g')$ denote the discrete sequence of the zeros of $g'$, repeated according to their multiplicity. Let $X$ be a subset of $\BMOA$, defined by
\[
X=\{G\in \BMOA:\; Z(g')\subseteq Z(G')\}
\]
For any function $G\in X$, we then have that $\frac{G'}{g'} \in \Hol(\mathbb D)$. It then follows that $X\subset \Lambda_a$, with $a=\frac{1}{p_1}-\frac{1}{p_2}$. To verify this, let $h\in H^{p_1}$. Then, the function $\psi=\frac{G'}{g'}\,h \in \HOl{\D}$, and
\[
T_g(\psi)(z)=\int_0^z g'(w)\,\frac{G'(w)}{g'(w)}\,h(w)\,dw\,=\,T_G(h)(z)\,,
\]
thus $T_g(\psi) \in H^{p_1}$, since $G\in \BMOA$. In different words, $\psi \in [T_g, H^{p_1}]$ which is equal to $[T_g, H^{p_2}]$ by our assumption. Therefore, we have that for every $h\in H^{p_1}$, $T_G(h) \in H^{p_2}$ which, in turn, means that
$T_G:\;H^{p_1}\to H^{p_2}$ is bounded. Then, in view of \eqref{E:AC}, $G$ must be in $\Lambda_a$.

However, $G \in \Lambda_a$ is equivalent to $G' \in K_{1-a}$, and applying Lemma \ref{L:missing}, we can find a function $\phi$, such that $\phi\cdot G' \in K_{1-a+\varepsilon} \setminus K_{1-a}$. Notice now, that for the function
\[
H_1(z)=\int_0^z \phi(w)\cdot G'(w)\,dw\,,
\]
we have that $H_1 \in X$, but $H'_1 \notin K_{1-a}$ or equivalently $H_1 \notin \Lambda_a$, which is a contradiction. 
\end{proof}

As a final remark, we note that finding a concrete example of a function $\phi \in [T_g,H^{p_1}]\setminus [T_g,H^{p_2}]$ seems challenging.

\section{Acknowledgments}
 The authors would like to express their gratitude to professors Petros Galanopoulos and Aristomenis Siskakis for useful discussions on this topic. We would also like to thank Nikolaos Chalmoukis for his insightful ideas that enriched this article.
\bibliographystyle{plain}

\bibliography{Literature}
\end{document}